\documentclass[twoside,12pt,a4paper]{article}
\usepackage{times}
\usepackage[all]{xy}
\usepackage{amssymb,amsbsy,amsmath,amsfonts,amssymb,amscd}
\usepackage{latexsym,euscript,exscale,epic,epsfig}
\begin{document}
\newtheorem{corollary}{Corollary} 
\newtheorem{result}{Theorem}
\newtheorem{little}{Lemma}
\newtheorem{proposition}{Proposition}
\newtheorem{definition}{Definition}
\newtheorem{conjecture}{Conjecture}
\newtheorem{question}{Question}
\newtheorem{remark}{Remark}
\newtheorem{def/thm}{Definition-Theorem}
\newtheorem{claim}{Claim}

\centerline{\huge Explicit Noether-Lefschetz for arbitrary threefolds}

\vskip 1cm

\centerline{\bf ANGELO FELICE LOPEZ and CATRIONA MACLEAN}

\begin{abstract} 
\noindent We study the Noether-Lefschetz locus of a very ample line bundle $L$ on an
arbitrary smooth threefold $Y$. Building on results of Green, Voisin
and Otwinowska, we give explicit bounds, depending only on the Castelnuovo-Mumford
regularity properties of $L$, on the codimension of the components of the
Noether-Lefschetz locus of $|L|$. 
\end{abstract}

\section{Introduction.}
It is well-known in algebraic geometry that the geometry of a given variety is influenced
by the geometry of its subvarieties. It less common, but not unusual, that a given
ambient variety forces to some extent the geometry of its subvarieties.
\newline
A particularly nice case of the latter is given by line
bundles, whose properties do very much influence the geometry. 
\newline If $Y$ is a smooth variety and $i : X \hookrightarrow Y$ is a smooth divisor, there is
then a natural restriction map
\[i^{\ast} : {\rm Pic}(Y) \rightarrow {\rm Pic}(X) \]
given by pull-back of line bundles.
\newline Now suppose that $X$ is very ample. By the Lefschetz theorem
$i^{\ast}$ is
injective if ${\rm dim} Y \geq 3$. On the other hand, it was already known to the Italian school
(Severi \cite{se}, Gherardelli \cite{ghe}), that $i^{\ast}$ is surjective when ${\rm dim} Y \geq
4$. Simple examples show that in the case where  ${\rm dim} Y = 3$ we cannot hope for surjectivity
unless a stronger restriction is considered. 
\newline 
For the case $Y = \mathbb{P}^3$, this is also a classical problem, first posed by Noether and
solved in the case of generic $X$ by Lefschetz who showed that
\newline\newline
{\bf Theorem (Noether-Lefschetz)}
{\it For $X$ a generic surface of degree $d \geq 4$ in $\mathbb{P}^3$ we have
${\rm Pic}(X)\cong \mathbb{Z}$.}
\newline\newline
Here and below by generic we mean outside a countable union of proper subvarieties.
\newline
Suppose now that a smooth threefold $Y$ and a line bundle $L$ on $Y$ are given.
We will say that a Noether-Lefschetz theorem holds for the pair $(Y, L)$, 
if 
\[i^*: {\rm Pic}(Y) \rightarrow {\rm Pic}(X)\] 
is a surjection for a generic smooth surface $X \subset Y$ such that $\mathcal{O}_Y(X)=L$.
\newline\newline
The following result of Moishezon (\cite{moi}, see also the argument given in Voisin
\cite[Thm.\ 15.33]{book}) establishes the exact conditions under which a Noether-Lefschetz 
theorem holds for $(Y,L)$.
\newline\newline
{\bf Theorem (Moishezon)}
{\it If $(Y,L)$ are such that $L$ is very ample and
\[h_{ev}^{0,2}(X, \mathbb{C}) \neq 0\]
for a generic smooth $X$ such that $\mathcal{O}_Y(X) = L$, then a Noether-Lefschetz theorem holds
for the pair $(Y, L)$.}
\newline\newline
Here, $h_{ev}^{0,2}$ denotes the evanescent $(2,0)$-cohomology of $X$: see below for a
precise definition.
\newline\newline
More precisely, we denote by $U(L)$ the open subset of $\mathbb{P}H^0(L)$ parameterizing smooth
surfaces in the same equivalence class as $L$. We further denote by ${\rm NL(L)}$ ({\it the
Noether-Lefschetz locus of $L$}) the subspace parameterizing surfaces $X$ equipped with line
bundles which are not produced by pull-back from $Y$. The above theorem then admits the 
following alternative formulation.
\newline\newline
{\bf Theorem (Moishezon)}
{\it If $(Y,L)$ are such that $L$ is very ample and
\[h_{ev}^{0,2}(X,\mathbb{C})\neq 0\]
for a generic smooth $X$ such that $\mathcal{O}_Y(X)=L$, then
the Noether-Lefschetz locus ${\rm NL(L)}$ is a countable union of proper algebraic 
subvarieties of $U(L)$.}
\newline\newline
These proper subvarieties will henceforth
be referred to as {\it components of the Noether-Lefschetz
locus}. 
\newline\newline
A Noether-Lefschetz theorem for a pair $(Y,L)$ essentially says that for a generic surface $X$
such that $\mathcal{O}_Y(X) = L$, the set of line bundles on $X$ is well-understood and as
simple as possible. A natural follow-up question is: how rare are surfaces with badly behaved
Picard groups? Or alternatively: how large can the components of the Noether-Lefschetz locus be
in comparison with $U(L)$? This leads us to attempt to prove what we will call {\it explicit
Noether-Lefschetz theorems}. An explicit  Noether-Lefschetz theorem (the terminology is due to
Green) says that the codimension of ${\rm NL(L)} \subset U(L)$ is bounded below by some number
$n_L$ depending non-trivially on the positivity of $L$. The first known example of these was the
following theorem, established independently by Voisin and Green, \cite{gd3}, \cite{vd3},
which gives an explicit Noether-Lefschetz theorem for $\mathbb{P}^3$.
\newline\newline
{\bf Theorem (Green, Voisin)}
{\it Let $Y = \mathbb{P}^3$ and $L = \mathcal{O}_{\mathbb{P}^3}(d)$. Let $\Sigma_L \subset U(L)$ be any
component of the Noether-Lefschetz locus. Then $ {\rm codim} \ \Sigma_L \geq d - 3$, with equality
being achieved only for the component of surfaces containing a line.}
\newline\newline
In this theorem we see also another of the reigning principles of the study of components of the
Noether-Lefschetz locus, namely that components of small codimension should parameterize
surfaces containing low-degree curves.
\newline\newline
Recently, the subject has been much advanced by the following result
of Ot\-wi\-now\-ska, 
(\cite{ania2}, see also \cite{ania} and \cite{ania1}) which implies an explicit Noether-Lefschetz theorem for analogues of
Noether-Lefschetz loci for highly divisible line bundles on varieties of arbitrary odd
dimension. (For ease of presentation, we give a weakened version of the result proved).
\newline\newline
{\bf Theorem (Otwinowska)}
{\it Let $Y$ be a projective variety of dimension $2n + 1$, let $\mathcal{O}_Y(1)$ be a very
ample line bundle on $Y$ and let $\Sigma_L \subset U({\mathcal{O}_Y(d))}$ be any component of
the Noether-Lefschetz locus. Let $X$ be a hypersurface contained in $\Sigma_L$. For $d$ large
enough, if 
\[{\rm codim} \ \Sigma_L \ \leq \frac{d^n}{n!}\]
then $X$ contains a $n$-dimensional linear space.}
\newline\newline
In fact, Otwinowska also gives a numerical criterion on $d$ and the codimension of $\Sigma_L$
under which $X$ necessarily contains a degree-$b$ n-dimensional subvariety.
\newline
We recall also the results of Joshi \cite{jo} and Ein-Lazarsfeld \cite[Prop. 3.4]{el}.
\newline\newline
The aim in this paper will be to shed light on the fact that it is the {\it Castelnuovo-Mumford
regularity properties} of a line bundle that insure that an explicit Noether-Lefschetz
theorem holds, independently on the divisibility properties.
\newline 
To state our first result we suppose that $Y$ is a smooth threefold and $H$ is a very ample
line bundle on $Y$. We define numbers $\alpha_Y$ and $\beta_Y$ as follows.
\begin{definition}
\label{alpha}
The integer $\alpha_Y$ is defined to be the minimal \underline{positive} integer such that $K_Y +
\alpha_Y H$ is very ample. The integer $\beta_Y$ is defined to be the minimal integer such that 
$(\beta_Y - \alpha_Y)H - K_Y$ is nef.
\end{definition}
We recall that, by the results of adjunction theory \cite{sv}, if $(Y, H) \neq
(\mathbb{P}^3, \mathcal{O}_{\mathbb{P}^3}(1))$, we have that $\alpha_Y \leq 4$ with equality if
and only if either $Y$ is a $\mathbb{P}^2$-bundle over a smooth curve and the restriction of $H$
to the fibers is $\mathcal{O}_{\mathbb{P}^2}(1)$ ({\it we will refer later to this case as a
linear $\mathbb{P}^2$-bundle}) or $(Y, H) = (Q, \mathcal{O}_Q(1))$ where $Q \subset \mathbb{P}^4$
is a smooth quadric hypersurface. On the other hand $\beta_Y \geq 1$ with equality if $Y$ is
subcanonical and nonpositive (that is if $K_Y = e H$ for some integer $e \leq 0$).
\newline 
We have
\begin{result}
\label{(-d)-regular}
Let $Y$ be a smooth threefold, $Y \neq \mathbb{P}^3$ and let $H$ be a very ample divisor on $Y$. Let
$L$ be a (-d)-regular line bundle with respect to $H$. We suppose that either $H^1(\Omega_Y^2 \otimes
L) = 0$ or $d \geq 3\beta_Y - 3\alpha_Y + 13$.  Let $\Sigma_L$ be any component of the
Noether-Lefschetz locus ${\rm NL(L)}$. The following bounds hold:
\begin{itemize}
\item[(i)] If $(Y, H)$ is not a linear $\mathbb{P}^2$-bundle then
\[{\rm codim} \ \Sigma_L \geq \begin{cases} d - 5 + \alpha_Y - 2\beta_Y & {\rm if} \ \beta_Y
\geq 2 \ {\rm and} \ d \geq \frac{\beta_Y^2 (\beta_Y + 5)}{2} \\ d - 6 + \alpha_Y &
{\rm if} \ \beta_Y = 1 \end{cases}. \]
\item[(ii)] If $(Y, H)$ is a linear $\mathbb{P}^2$-bundle then
\[{\rm codim} \ \Sigma_L \geq \begin{cases} d - 2 - 2\beta_Y & {\rm if} \ \beta_Y \geq 2 \ 
{\rm and} \ d \geq \frac{\beta_Y^2(\beta_Y + 5)}{2} \\ d - 3 & {\rm if} \
\beta_Y = 1 \end{cases}. \]
\end{itemize}
\end{result}
We can do a little bit better in the case of the Noether-Lefschetz locus of adjoint line bundles.
\newline\newline
We now define numbers $a_Y$ and $b_Y$ as follows.
\begin{definition}
\label{a}
The integer $a_Y$ is defined to be the minimal integer such that $K_Y + a_YH$ is very ample. 
The integer $b_Y$ is defined to be the minimal integer such that $(b_Y - a_Y)H - K_Y$ is nef.
\end{definition}
As above, if $(Y, H) \neq (\mathbb{P}^3, \mathcal{O}_{\mathbb{P}^3}(1))$, we have that $a_Y \leq
4$ with equality if and  only if either $(Y, H)$ is a linear $\mathbb{P}^2$-bundle or $(Y, H) =
(Q, \mathcal{O}_Q(1))$ and again $b_Y \geq 1$ with equality if $Y$ is subcanonical.\newline
\begin{result} 
\label{pluriadjoint} 
Let $Y$ be a smooth threefold, $Y \neq \mathbb{P}^3$ and let $H$ be a very ample divisor on $Y$. Let
\[L = K_Y + dH + A,\] 
where $A$ is numerically effective.
We suppose that either $H^1(\Omega_Y^2 \otimes L) = 0$ or $d \geq 2b_Y - 2a_Y + 13$.   
Let $\Sigma_L$ be any component of the Noether-Lefschetz locus ${\rm NL(L)}$. The following
bounds hold:
\begin{itemize}
\item[(i)] If $(Y, H)$ is not a linear $\mathbb{P}^2$-bundle then
\[{\rm codim} \ \Sigma_L \geq \begin{cases} d - 5 - b_Y & {\rm if} \ b_Y \geq 2 \ {\rm and} \ 
d \geq \frac{b_Y (b_Y^2 + 7 b_Y -6)}{2} \\ d - 5 & {\rm if} \ b_Y = 1 \end{cases}. \]
\item[(ii)] If $(Y, H)$ is a linear $\mathbb{P}^2$-bundle then
\[{\rm codim} \ \Sigma_L \geq \begin{cases} d - 6 - b_Y & {\rm if} \ b_Y \geq 2 \ {\rm and} \ 
d \geq \frac{b_Y(b_Y-1)(b_Y + 8)}{2} \\ d - 6 & {\rm if} \ b_Y = 1 \end{cases}. \]
\end{itemize}
\end{result}
We also note the following application that generalises \cite{an} (see also \cite{ba}).
\begin{corollary}
\label{blow-up}
Let $Y$ be a smooth threefold such that $Y \neq \mathbb{P}^3$ and ${\rm Pic}(Y) \cong \mathbb{Z}
H$ where $H$ is a very ample line bundle and let $K_Y = eH$. We suppose that either $H^1(\Omega_Y^2
(d)) = 0$ or $d \geq 3e + 13$. Let $P_1, \ldots, P_k$ be $k$ general points in $Y$ and $\pi :
\widetilde{Y} \to Y$ be the blow-up of $Y$ at these points with exceptional divisors $E_1, \ldots,
E_k$.
\newline
If $d \geq 7 + e$ then
\[d \pi^{\ast}(H) - E_1 - \ldots - E_k \ \ \ {\rm \ is \ ample \ on} \ \widetilde{Y}
\ \Leftrightarrow \ \ d^3 H^3 > k .\]
\end{corollary}
We outline our approach to the study of the Noether-Lefschetz locus.
\newline
In section 2, we will give the standard expression of this problem in terms of variation of
Hodge structure of $X$. We will then recall the classical results of Griffiths, Carlson et. al.
which allow us to express variation of Hodge structure of $X$ in terms of multiplication of
sections of line bundles on $X$. 
\newline\newline
We define $\sigma$ to be the section of $L$ defining $X$. The tangent space of a component of 
the Noether-Lefschetz locus is naturally a subspace of $H^0(L)/\langle \sigma\rangle$, and we
will denote its preimage in $H^0(L)$ by $T$. If we suppose that $H^1(\Omega_Y^2 \otimes L) =
0$, then $T$ has the following property: The natural multiplication map
\begin{equation}
\label{main} 
T \otimes H^0(K_Y \otimes L) \rightarrow H^0(K_Y \otimes L^2)
\end{equation}
is not surjective.\newline
A full proof of this fact is given in section 3.
\newline
In section 3, we also explain Green's methods for proving the explicit Noether-Lefschetz theorem
for $\mathbb{P}^3$ using Koszul cohomology to prove that equation (\ref{main}) cannot be
satisfied if $T$ is too large. Green's method does not immediately
apply to our case, since it requires $T$ to be base-point free--- which is only guaranteed if the
tangent bundle of $Y$ is globally generated, hence only for a few threefolds. However, we show in
section 4 that there exists $W \subset H^0(K_Y \otimes L(3))$ such that $W$ is base-point free and 
\[\{T \otimes H^0(K_Y \otimes L) \} \oplus \{W \otimes H^0(L(-3)) \} \rightarrow H^0(K_Y \otimes L^2) \]
is not surjective. Results of Ein and Lazarsfeld \cite{el} 
then imply a lower bound on the codimension of 
\[\{T \otimes H^0(K_Y(3)) \} \oplus W \subset H^0(K_Y \otimes L(3)) \]
and more particularly on the codimension of
\[T \otimes H^0(K_Y(3)) \subset H^0(K_Y \otimes L(3)).\]
In introducing $W$, we get around the base-point free problems, but introduce others. In
particular, we now need a method for extracting a lower bound on ${\rm codim} \ T$ from a lower
bound for ${\rm codim} \ (T \otimes H^0(K_Y(3)))$. When $Y = \mathbb{P}^3$, this is a simple 
application of a classical inequality in commutative algebra due to Macaulay and Gotzmann.
In section 5 we extend the Macaulay-Gotzmann inequality to sections of any Castelnuovo-Mumford
regular sheaf. In section 6, we pull all of the above together to prove the theorem.
\section{Preliminaries.}
In this section we resume the classical results of Griffiths, Carlson
et. al. on which our work will be based. We will show how a component $\Sigma_L$ of the
Noether-Lefschetz locus ${\rm NL(L)}$ can be locally expressed as the zeros of a certain section
of a vector bundle over $U(L)$. We will then use this expression--- together with the work of
Griffiths from the 60s, relating variation of Hodge structure with deformations of $X$
to multiplication of sections of line bundles on $X$--- to relate the
codimension of $\Sigma_L$ to cohomological questions on $X$.  
\subsection{${\rm NL}$ expressed as the zero locus of a vector bundle section.} 
We note first that by the Lefschetz theorem the map ${\rm Pic}_0(Y)\rightarrow {\rm Pic}_0(X)$ is
necessarily  an isomorphism. It follows that the map 
\[i^*:{\rm Pic}(Y)\rightarrow {\rm Pic}(X)\] 
fails to be surjective if and only if the $(1,1)$ integral evanescent cohomology is non-trivial:
\[H^{1,1}_{\rm ev}(X, \mathbb{Z}) \neq 0.\] 
(We recall that the subspace $H^{1,1}_{\rm ev}(X, \mathbb{C}) \subset H^{1,1}(X, \mathbb{C})$ is
defined by
\[\gamma \in H^{1,1}_{\rm ev}(X, \mathbb{C}) \Leftrightarrow \langle i^*\beta, \gamma\rangle
= 0 \mbox{ for all } \beta \in H^2(Y,\mathbb{C}).) \]
In particular, we can therefore define ${\rm NL(L)}$ as follows
\[X \in {\rm NL(L)} \Leftrightarrow H^{1,1}_{\rm ev}(X,\mathbb{Z}) \neq 0.\] 
This is the definition of ${\rm NL(L)}$ which we will use henceforth, since it is much more
manageable. In particular, it is this description which will allow us to write any component of
${\rm NL(L)}$ as the zero locus of a special section of a vector bundle.
\newline\newline
Henceforth, we will assume that $X$ is contained in ${\rm NL(L)}$ and $\gamma$ will be a
non-trivial element of $H^{1,1}_{\rm ev}(X, \mathbb{Z})$. The point in $U(L)$
corresponding to $X$ will be denoted by $0$. We will now define what we mean by the {\it
Noether-Lefschetz locus associated to $\gamma$}, which we denote by $NL(\gamma)$. Since we will
be interested in the local geometry of ${\rm NL(L)}$, we fix for simplicity a contractible
neighbourhood of $0$, $O$. Henceforth, all our calculations  will be made over $O$. We form a
vector bundle $\mathcal{H}_{\rm ev}^2$ over $O$, defined by
\[\mathcal{H}_{\rm ev}^2(u)= H^2_{\rm ev}(X_u,\mathbb{C}).\]
The vector bundle contains holomorphic sub-bundles $\mathcal{F}^i(\mathcal{H}_{\rm ev}^2)$ given
by 
\[\mathcal{F}^i(\mathcal{H}_{\rm ev}^2)(u)=F^i(H_{\rm ev}^2(X_u,\mathbb{C})).\]
We define bundles $\mathcal{H}_{\rm ev}^{i, 2-i}$ by
\[\mathcal{H}_{\rm ev}^{i, 2-i}= \mathcal{F}^i(\mathcal{H}_{\rm ev}^2)/
\mathcal{F}^{i+1}(\mathcal{H}_{\rm ev}^2).\] 
(The fibre of $\mathcal{H}_{\rm ev}^{i, 2-i}$ at the point $u$ is isomorphic to $H_{\rm ev}^{i,
2-i}(X_u)$:  however, $\mathcal{H}_{\rm ev}^{i,2-i}$ does not embed naturally into  
$\mathcal{H}_{\rm ev}^2$ as a holomorphic sub-bundle.) The bundle $\mathcal{H}_{\rm ev}^2$ is
equipped with a natural flat connexion, the Gauss-Manin connexion, which we denote by $\nabla$.
We now define $\overline{\gamma}$ to be the section of $\mathcal{H}_{\rm ev}^2$ produced by flat
transport of $\gamma$.
\newline\newline
We define $\overline{\gamma}^{0,2}$, a section of $\mathcal{H}_{\rm ev}^{0,2}$, to be the 
image of $\overline{\gamma}$ under the projection
\[\pi:\mathcal{H}_{\rm ev}^2\rightarrow \mathcal{H}_{\rm ev}^{0,2}.\] 
We are now in a position to define ${\rm NL}(\gamma)$.
\begin{definition}
The Noether-Lefschetz locus associated to $\gamma$, ${\rm NL}(\gamma)$, 
is given by
\[{\rm NL}(\gamma)= {\rm zero}(\overline{\gamma}^{0,2}).\]
\end{definition}
Informally, ${\rm NL}(\gamma)$ parameterizes the small deformations of $X$ on which $\gamma$
remains of Hodge type $(1,1)$. Any component of ${\rm NL(L)}$ is locally equal to ${\rm
NL}(\gamma)$ for some $\gamma$. 
\newline\newline
The tangent space $T{\rm NL}(\gamma)$ at $X$ is a subspace of $H^0(L)/\langle \sigma \rangle$,
where $\sigma$ is the section of $L$ defining $X$. We will  denote its preimage in $H^0(L)$ by
$T$. 
\subsection{IVHS and residue maps.}
We will now explain the classical work of Griffiths which makes the section
$\overline{\gamma}^{0,2}$ particularly manageable.
\newline\newline
Let $\mathcal{H}_{\rm ev}^2$ be as above. For the purposes of this section we will consider the
holomorphic subvector bundle $\mathcal{F}_{\rm ev}^i$ to be a holomorphic map
\[\mathcal{F}_{\rm ev}^i: O \rightarrow {\rm Grass}(f_i, \mathcal{H}_{\rm ev}^2) \]
where $f_i$ is the dimension of $F^iH_{\rm ev}^2(X,\mathbb{C})$. The Gauss-Manin connexion
gives us a canonical isomorphism
\[\mathcal{H}_{\rm ev}^2\cong H_{\rm ev}^2(X,\mathbb{C})\times O\]
from which we deduce a canonical isomorphism
\[ {\rm Grass}(f_i, \mathcal{H}_{\rm ev}^2)\cong O\times  
{\rm Grass}(f_i, H_{\rm ev}^2(X,\mathbb{C})).\]
In particular,  $\mathcal{F}_{\rm ev}^i$ is now expressed as a map from $O$ to the constant
space  
\newline
${\rm Grass}(f_i, H_{\rm ev}^2(X,\mathbb{C}))$, and as such can be derived.  We obtain a
derivation map, which we denote by ${\rm IVHS}$ (for Infinitesimal Variation of Hodge Structure)
\[{\rm IVHS}^i: TO\rightarrow {\rm Hom}(F^i(H_{\rm ev}^2),
H_{\rm ev}^2/F^i(H_{\rm ev}^2)).\]
Griffiths proved the following result in \cite{trans}. 
\newline\newline
{\bf Theorem (Griffiths' Transversality)}
{\it The image of ${\rm IVHS}^i$ is contained in
\newline 
${\rm Hom}(H_{\rm ev}^{i, 2-i}, H_{\rm ev}^{i-1, 3-i})$.}
\newline\newline
The importance of this work for our purposes is the following lemma.
\begin{little}
For any $v \in TO$, we have that
\[d_v(\overline{\gamma}^{0,2}) = -{\rm IVHS}^1(v)(\gamma) \]
\end{little}
{\bf Proof.} 
The isomorphism 
\[f:T_W{\rm Grass(n, V)} \cong {\rm Hom}(W, V/W)\] 
is given by
\[f(v): w\rightarrow \frac{\partial}{\partial v}(\tilde{w})_{|_{V/W}} \]
where $w \in W$ and $\tilde{w}$ is any local section of the tautological
bundle over the Grassmannian such that $\tilde{w}_W=w$.
\newline
In the case in hand,  we choose a lifting of $\overline{\gamma}^{0,2}$ to a section of 
$\mathcal{H}_{\rm ev}^2$, which we denote by $\overline{\gamma}^{0,2}_{\rm lift}$.
By definition of $\overline{\gamma}^{0,2}$, we then have that 
\[\overline{\gamma}-{\overline{\gamma}}^{0,2}_{\rm lift} \in \mathcal{F}^1
(\mathcal{H}_{\rm ev}^2)\]
and it follows that
\[{\rm IVHS}^1(v)(\gamma)= \frac{\partial}{\partial v}
(\overline{\gamma}-\overline{\gamma}^{0,2}_{\rm lift})_{|_{H_{\rm ev}^{0,2}}} \]
and now, since by definition $\overline{\gamma}$ is constant,
\[\hskip 3cm {\rm IVHS}^1(v)(\gamma) = -d_v(\overline{\gamma}_{\rm lift}^{0,2})_{|_{H_{\rm
ev}^{0,2}}} = -d_v(\overline{\gamma}^{0,2}). \hskip 2.8cm \square\]  
\newline\newline
We will also need the work of Carlson and Griffiths relating
the residue maps to Hodge structure of varieties (\cite{cg}). Suppose given, for $i = 1, 2$,
a section 
\[s \in H^0(K_Y \otimes L^i).\] 
This can be thought of as a holomorphic 3-form on $Y$ with a pole of order $i$ along $X$, 
and as such defines a cohomology class in $H^3(Y\setminus X, \mathbb{C})$. The group 
$H^3(Y\setminus X, \mathbb{C})$ maps to $H_{\rm ev}^2(X,\mathbb{C})$ via
residue, and hence there is an induced residue map
\[{\rm res_i}:  H^0(K_Y \otimes L^i)\rightarrow H_{\rm ev}^2(X,\mathbb{C}). \]
The relevance of this map to variation of Hodge structure comes from
the following theorem, which is proved by Griffiths in \cite{res}.
\newline\newline
{\bf Theorem}
{\it The image of ${\rm res_i}$ is contained in $F^{3-i}(H_{\rm ev}^2)$.}
\newline\newline
Henceforth, we will denote by $\pi_i$ the induced projection map
\[\pi_i: H^0(K_Y \otimes L^i)\rightarrow H_{\rm ev}^{3-i, i-1}
(X,\mathbb{C}).\]
In this representation, the map ${\rm IVHS}^{3-i}$ has a particularly nice form (\cite{cg}, page 70). 
\newline\newline
{\bf Theorem (multiplication)}
{\it Consider $v \in TO$. Let $\tilde{v}$ be a lifting of 
$v$ to $H^0(L)$. Then for any $P \in H^0(K_X\otimes L^i)$, we have that
\[{\rm IVHS}^{3-i}(v)(\pi_i (P)) = \pi_{i+1}(\tilde{v} \otimes P)\]
up to multiplication by some nonzero constant.}
\newline\newline
The only fly in the ointment is that in general we cannot be sure that 
the map $\pi_i$ is surjective onto $H_{\rm ev}^{3-i, i-1}(X,\mathbb{C})$. It is
precisely for this reason that we will be obliged to suppose that 
$H^1(\Omega_Y^2 \otimes L) = 0$.
\newline
The following lemma will be crucial.
\begin{little}
\label{orthog}
Consider $\gamma \in H_{\rm ev}^{1,1}(X)$ and $\omega \in H_{\rm ev}^{2,0}(X)$. 
For any vector $v \in TO$ we have
\[\langle {\rm IVHS}^1(v)(\gamma), \omega\rangle +
\langle \gamma,  {\rm IVHS}^2(v) (\omega)\rangle=0\]
\end{little}
{\bf Proof.} We note that
\[d_v(\langle \overline{\gamma}, \overline{\omega}\rangle)=0.\]
We note that we can write
\[\overline{\gamma}= \overline{\gamma}^1+\overline{\gamma}^2\]
where $\overline{\gamma^1} \in \mathcal{F}_{\rm ev}^1$ and 
$\overline{\gamma}^2(0)=0$. Similarly, we can write 
\[\overline{\omega}= \overline{\omega}^1+\overline{\omega}^2\]
where $\overline{\omega^1} \in \mathcal{F}_{\rm ev}^2$ and 
$\overline{\omega}^2(0)=0$. We note that for Hodge theoretic reasons
\[<\overline{\omega}^1, \overline{\gamma}^1>=0\]
and hence
\[d_v(\langle \overline{\gamma}, \overline{\omega}\rangle)
= \langle d_v(\overline{\gamma}^2),\omega\rangle
+ \langle \gamma, d_v(\overline{\omega}^2)\rangle.\]
Here, of course, it makes sense to talk about
$d_v(\overline{\omega}^2)$
and $d_v(\overline{\gamma}^2)$ only because $\overline{\omega}^2(0)=0$
and $\overline{\gamma}^2(0)=0$. 
Since $\langle \mathcal{F}^1,\mathcal{F}^2\rangle=0$, we have that
\[\langle d_v(\overline{\gamma}^2), \omega\rangle 
=\langle d_v(\overline{\gamma}^2)^{0,2}, \omega \rangle
=\langle -{\rm IVHS}^1(v)(\gamma),\omega\rangle\]
and similarly
\[\langle \gamma, d_v(\overline{\omega}^2)\rangle
=\langle \gamma, (d_v\overline{\omega}^2)^{1,1}\rangle
=\langle \gamma, -{\rm IVHS}^2(v)(\omega)\rangle.\]
So it follows immediately from 
\[d_v(\langle \overline{\gamma}, \overline{\omega}\rangle)=0\]
that 
\[\hskip 3cm \langle {\rm IVHS}^1(v)(\gamma), \omega\rangle +
\langle \gamma,  {\rm IVHS}^2(v) (\omega)\rangle=0. \hskip 2.8cm \square\] 
\section{Strategy and overview.}
The basic idea of this proof is that used by Green in \cite{gd3}. 
We summarise his proof and explain why it cannot be immediately applied to the situation in 
hand.
\newline\newline
First some notation. Given any pair of coherent sheaves on $X$, $L$ and $M$ we denote by 
$\mu_{L,M}$ the multiplication map
\[\mu_{L,M}:H^0(L)\otimes H^0(M)\rightarrow H^0(L\otimes M).\]
Where there is no risk of confusion, we will write $\mu$ for $\mu_{L,M}$.
\newline\newline
The starting point of Green's work is the following lemma.
\begin{little}
\label{green}
Suppose that $T \subset H^0(\mathcal{O}_{\mathbb{P}^3}(d))$ is the preimage of 
$TNL(\gamma)$. Then the inclusion
\[\mu(T \otimes H^0(\mathbb{P}^3, \mathcal{O}_{\mathbb{P}^3}(d - 4)))
\subset H^0(\mathbb{P}^3, \mathcal{O}_{\mathbb{P}^3}(2 d - 4))\]
is a \underline{strict} inclusion.
\end{little}
{\bf Proof.}
In the case of $Y = \mathbb{P}^3$, we have that $\pi_i: H^0(K_Y \otimes L^i) \rightarrow H_{\rm
ev}^{3-i, i-1}(X)$  is a surjection. (See, for example, \cite[proof of Thm.\ 18.5, page
420]{book}). By Lemma \ref{orthog}, if $v \in TNL(\gamma)$ and 
$P \in H^0(\mathbb{P}^3, \mathcal{O}_{\mathbb{P}^3}(d-4))$ then
\[\langle \gamma, {\rm IVHS}^2(v)(\pi_1(P)) \rangle= -\langle   
{\rm IVHS}^1(v)(\gamma), \pi_1(P) \rangle = 0\]
from which we conclude that 
\[{\rm IVHS}^2(v)(\pi_1(P)) \in \gamma^\perp, \]
where $\gamma^\perp$ is the orthogonal to $\gamma$, and in particular
is a proper subspace.
\newline 
By the multiplication theorem it follows that
\[\pi_2(\mu(\tilde{v} \otimes P)) \in \gamma^\perp \] 
or alternatively
\[\mu(\tilde{v} \otimes P) \in \pi_2^{-1}(\gamma^\perp).\] Since $\pi_2$ is 
surjective, $\pi_2^{-1}(\gamma^\perp)$ is a proper subspace. $\hfill\square$ 
\newline\newline
Green then proves the following theorem via the vanishing of certain Koszul cohomology groups.
\newline\newline
{\bf Theorem (Green)}
{\it Let $T \subset H^0(\mathcal{O}_{\mathbb{P}^r}(d))$ be a base-point free linear system of 
codimension $c$. Then the Koszul complex 
\[\bigwedge^{p+1} T \otimes H^0(\mathcal{O}_{\mathbb{P}^r}(k-d)) \rightarrow 
\bigwedge^p T \otimes H^0(\mathcal{O}_{\mathbb{P}^r}(k)) \rightarrow 
\bigwedge^{p-1} T \otimes H^0(\mathcal{O}_{\mathbb{P}^r}(k+d)) \]
is exact in the middle provided that $k\geq p+d+c$.}
\newline\newline
In the case in hand, on setting $r = 3, p=0$ and $k = 2d - 4$ we see that the multiplication map
\[T \otimes H^0(\mathbb{P}^3,\mathcal{O}_{\mathbb{P}^3}(d - 4)) \rightarrow
H^0(\mathbb{P}^3, \mathcal{O}_{\mathbb{P}^3}(2d - 4))\] 
is surjective if $2d - 4 \geq d + c$. But we have already observed that this multiplication 
map is necessarily non-surjective, from which it follows that $c \geq d - 3$.
\newline\newline
In Lemma \ref{omega2} below we will see that, provided $H^1(\Omega_Y^2 \otimes L) = 0$, 
it is still true that the multiplication map
$T \otimes H^0(K_Y \otimes L) \rightarrow H^0(K_Y \otimes L^2)$ is non-surjective.
One might therefore reasonably entertain the hope of adapting Green's methods to arbitrary 
varieties. The difficulty is that in order to apply Green's result, $T$ must be base-point 
free. This was immediate when $Y = \mathbb{P}^3$, since, if $X$ was given by $F \in 
H^0(\mathcal{O}_{\mathbb{P}^3}(d))$, $T$ then automatically  contained
$H^0(\mathcal{O}_{\mathbb{P}^3}(1)) \times \langle \frac{\partial F}{\partial X_i} \rangle$. However if
$T_Y$ is not globally  generated, there is no reason why this should hold in general. The rest of this
paper will be concerned with finding ways around this difficulty.
\begin{little}
\label{omega2}
Let $L$ be very ample and such that 
\[H^1(\Omega_Y^2 \otimes L) = 0.\]
Let $T \subset H^0(L)$ be the preimage in $H^0(L)$ of the tangent space to
${\rm NL}(\gamma)$. Then 
\[\mu(T\otimes H^0(K_Y \otimes L)) \subset H^0(K_Y \otimes L^2)\]
is a strict inclusion. 
\end{little}
{\bf Proof.}
We note that by the argument given in the proof of Lemma \ref{green},
\[\pi_2(\mu(T\otimes H^0(K_Y \otimes L))) \neq H^{1,1}_{\rm ev}(X,\mathbb{C}).\]
Now it just remains to observe that, by \cite[proof of Thm.\ 18.5, page 420]{book}, 
\[\pi_2: H^0(K_Y \otimes L^2)\rightarrow H^{1,1}_{\rm ev}(X,\mathbb{C})\] 
is a surjection, since
\[\hskip 4.4cm H^1(\Omega_Y^2(X))=0. \hskip 5.8cm \square\]  
\newline
So, we would now like to apply Green's argument; unfortunately, $T$
may have base points. Our strategy for getting around this problem will be
as follows.
\begin{enumerate}
\item First of all, we will construct $W \subset H^0(K_Y \otimes L(3))$ 
with the following 
good properties.
\begin{enumerate}
\item $W$ is base-point free,
\item $\pi_2(\mu(W\otimes H^0(L(-3))))=0$.
\end{enumerate}
\item The result proved by Ein and Lazarsfeld in \cite{el}
then gives us a lower bound on the codimension of $\mu(T \otimes H^0(K_Y(3)))$ in
$H^0(K_Y \otimes L(3))$.
\item We will then extract from the lower bound on ${\rm codim} \ \mu(T \otimes H^0(K_Y(3)))$
a  lower bound on the codimension of $T$ in $H^0(L)$. 
\end{enumerate}  
\section{Constructing $W$.}
We henceforth let $Y$ be a smooth threefold, $Y \neq \mathbb{P}^3$ and $H$ be a very ample divisor on
$Y$.
\begin{proposition}
\label{exofW}
There is a subspace $W \subset H^0(K_Y \otimes L(3))$ such that 
\begin{enumerate}
\item  The map $\pi_2 \circ \mu: W \otimes H^0(L(-3)) \rightarrow 
H^{1,1}_{\rm ev}(X, \mathbb{C})$ is identically zero. 
\item $W$ is base-point free.
\end{enumerate}
\end{proposition}
{\bf Proof.}
We denote the image of \[\mu:W\otimes H^0(L(-3))\rightarrow 
H^0(K_Y \otimes L^2)\] by $\langle W\rangle$. Consider the map
\[d:  H^0(\Omega^2_Y\otimes L)\rightarrow H^0(K_Y \otimes L^2)\] 
which sends a two-form on $Y$ with a simple pole along $X$ to  its derivation. We note that for
any $\omega \in H^0(\Omega^2_Y\otimes L)$ we have that
\[d\omega \in {\rm Ker}({\rm res}_2),\] because $d\omega$, being exact,  
defines a null cohomology class on $Y \setminus X$.
\newline
The space $W$ will be chosen in such a way that
\[\langle W \rangle_{|_X} \subset {\rm Im}(d)_{|_X}.\]
The map $d$ is difficult to deal with because it is not a map of
$\mathcal{O}_Y$-modules: the value of $d \omega$ at a point $x$ is
not determined by the value of $\omega$ at $x$.  In particular, it is not possible to
form a tensor product map
\[d\otimes (L^{-1}(3)) : H^0( \Omega^2_Y(3)) \rightarrow  H^0(K_Y \otimes L(3)).\]
Our first step will be to show that, even if $d$ does not come from an underlying  map
of $\mathcal{O}_Y$-modules, the restriction
\[d_X: H^0 (\Omega^2_Y \otimes L)\rightarrow H^0(K_X \otimes L_{|_X}).\]
does.
\begin{little}
\label{d_X}
Let the map 
\[r: \Omega^2_Y \otimes L \rightarrow K_X \otimes L\] 
be given by tensoring with $L$ the pull-back 
\[i^* : \Omega^2_Y \rightarrow \Omega_X^2 (\cong K_X).\] 
Then we have that 
\[d_X = -H^0(r).\]
\end{little}
{\bf Proof.}
We calculate in analytic complex co-ordinates near a point $p \in X$. 
Let $f$ be a function defining $X$ in a  neighbourhood of $p$ and
let $x,y$ be co-ordinates chosen in such a way that $(f,x,y)$ form a 
system of co-ordinates for $Y$ close to $p$. If $\nu \in H^0(\Omega_Y^2 \otimes L)$, then in a
neighbourhood of $p$ we can write
\[ \nu = \frac{f_1 dx \wedge dy + f_2 dx \wedge df+ f_3 dy \wedge df}{f}\]
where $f_1, f_2, f_3$ are holomorphic functions on a neighbourhood of $p$.
\newline
Differentiating and restricting to $X$, we get that
\[ d\nu_{|_X} = \frac{- f_1 dx \wedge dy \wedge df}{f^2}. \]
As an element of $H^0((K_Y \otimes L) \otimes L)$, this is represented by
\[\frac{-f_1 dx\wedge dy\wedge df}{f} \otimes {1/f}.\] 
Under the canonical isomorphism $(K_Y \otimes L)_{|_X} \rightarrow K_X$, we have that
\[\frac{- f_1 dx \wedge dy \wedge df}{f} \rightarrow - f_1 dx \wedge dy.\] 
Hence, under the canonical isomorphism 
\[(K_Y \otimes L^2)_{|_X} \rightarrow K_X\otimes L_{|_X},\] 
we have that
\[(d\nu)_{|_X} \rightarrow \frac{- f_1 dx \wedge dy}{f} = - r(\nu).\]
This concludes the proof of Lemma \ref{d_X}. $\hfill\square$
\newline\newline
We now proceed with the proof of Proposition \ref{exofW}.
\newline
The map $d_X$, which \underline{is} a map of $\mathcal{O}_Y$-modules, has the advantage 
that we can form tensor products. We consider the map induced by tensor product with 
$L^{-1}(3)$
\[d_X^{L^{-1}(3)} : H^0(\Omega^2_Y(3)) \rightarrow H^0(K_X(3)).\]
We define $W$ by
\[W = \{w \in H^0(K_Y \otimes L(3)) : w_{|_X} \in {\rm Im}
(d_X^{L^{-1}(3)})\}.\]
We will prove first that
\begin{little} 
\label{Pw}
For any $w \in W$ and $P \in H^0(L(-3))$, we have that
\[\pi_2(\mu(P \otimes w)) = 0.\]
\end{little}
{\bf Proof.} 
Since $w \in W$ there exists $s \in H^0(\Omega_Y^2(3))$ such that
\[w_{|_X} = d_X^{L^{-1}(3)}s\]
and hence
\[(Pw)_{|_X} = d_X(Ps) = d(Ps)_{|_X}.\]
From this it follows that there exists $s' \in H^0(K_Y \otimes L)$ such that
\[Pw = d(Ps) + \sigma s'.\]
We observed above that $\pi_2(d(Ps))=0$. We note that
\[{\rm res}_2(\sigma s') = {\rm res}_1(s') \] 
and hence 
\[{\rm res}_2(\sigma s') \in F^2 H_{\rm ev}^2(X, \mathbb{C}),\] 
from which it follows that $\pi_2(\sigma s') = 0.$ Whence 
\[\pi_2(Pw) = 0.\]
This concludes the proof of Lemma \ref{Pw}. $\hfill\square$
\newline\newline
To conclude the proof of Proposition \ref{exofW} it remains only to show that $W$ is base-point
free. Since $Y \neq \mathbb{P}^3$ we have (\cite{ein}) that $K_Y(3)$ is globally generated. Also 
\[\mu(\mathbb{C} \sigma \otimes H^0(K_Y(3))) \subset W \]
therefore the only possible base points of $W$ are the points of $X$.
Consider an arbitrary point $p \in X$. Now if $\mathbb{P}^N = \mathbb{P}H^0(Y, H)$ we have that
$\Omega_Y^2 (3)$ is globally generated since $\Omega_{\mathbb{P}^N}^2(3)$ is such and there is a
surjection $\Omega_{\mathbb{P}^N}^2 (3) \twoheadrightarrow \Omega_Y^2 (3)$. Whence there exists
a section
\[s \in H^0(\Omega_Y^2(3)) \] 
such that $d_X^{L^{-1}(3)}(s)(p) \neq 0$. From the short exact sequence
\[0 \rightarrow K_Y(3) \rightarrow K_Y \otimes L(3) \rightarrow K_X(3) \rightarrow 0\] 
and Kodaira vanishing we see that there exists 
\[w \in H^0(K_Y \otimes L(3))\] 
such that $w_{|_X} = d_X^{L^{-1}(3)}(s)$. It follows that $w \in W$, and 
\[w(p) = d_X^{L^{-1}(3)}(s)(p) \neq 0.\] 
Hence $p$ is not a base-point of $W$.  
This completes the proof of Proposition \ref{exofW}. $\hfill\square$
\newline\newline
To get lower bounds on the codimension we will apply the following result of Ein and 
Lazarsfeld, \cite[Prop.\ 3.1]{el}.
\newline\newline
{\bf Theorem (Ein, Lazarsfeld)}
{\it Let $H$ be a very ample line bundle and $B, C$ be nef line bundles on a smooth complex
projective $n$-fold $Z$. We set
\[F_f = K_Z + f H + B \mbox{ and } G_e = K_Z + e H + C.\]
Let $V \subset H^0(Z, F_f)$ be a base-point free subspace of codimension $c$ and consider the
Koszul-type complex
\[\bigwedge^{p+1} V \otimes H^0(G_e) \rightarrow \bigwedge^p V \otimes 
H^0(F_f + G_e) \rightarrow \bigwedge^{p-1} V \otimes H^0(2F_f + G_e).\] 
If $(Z, H, B) \neq (\mathbb{P}^n, \mathcal{O}_{\mathbb{P}^n}(1), \mathcal{O}_{\mathbb{P}^n})$, $f \geq
n+1$ and $e \geq n+p+c$, then this complex is exact in the middle.}
\newline\newline
In order to apply this to our situation, we set $p = 0$, and, in case $L = K_Y + dH + A$ we
choose $f = d, e = d - 3$, $B = A + K_Y + 3H$ (note that $B$ is nef since $K_Y + 3H$ is
globally generated) and $C = A$. In the case $L$ (-d)-regular we have $L = M(d)$ for a
Castelnuovo-Mumford regular line bundle $M$ and we choose $f = d + 3, e = d - 3 + \alpha_Y -
\beta_Y$, $B = M$ and $C = M + (\beta_Y - \alpha_Y)H - K_Y$, so that $B$ is nef since $M$ is
globally generated and also $C$ is nef by definition of $\alpha_Y$ and $\beta_Y$. We then have
that
\[F_f = K_Y \otimes L(3) \mbox{ and } G_e = L(-3) \]
and the theorem in this particular case says that:
\begin{proposition}
\label{el}
Suppose that $d \geq 4$ and $Y \neq \mathbb{P}^3$. Let $V$ be a base-point free linear system
in $H^0(K_Y \otimes L(3))$  with the property that
\[\mu(V \otimes H^0(L(-3))) \subset H^0(K_Y \otimes L^2) \]
is a strict inclusion. Then the codimension $c$ of $V$
satisfies the inequality
\[c \geq \begin{cases} d - 5 + \alpha_Y - \beta_Y & {\rm if \ L \ is \ (-d)-regular} \\ d - 5 &
{\rm if} \ L = K_Y + dH + A \end{cases}.\]
\end{proposition}
In general, pulling together the results of sections 3 and 4, we have the following bound.
\begin{proposition} 
\label{firstlower}
Suppose that $Y \neq \mathbb{P}^3$ and $H^1(\Omega_Y^2 \otimes L) = 0$. Then
the codimension of the image of
\[\mu: T \otimes H^0(K_Y(3)) \rightarrow H^0(K_Y \otimes L(3))\]
is at least $d - 5 + \alpha_Y - \beta_Y$ if $L$ is (-d)-regular or at least $d - 5$ if $L = K_Y
+ dH + A$.
\end{proposition}
{\bf Proof.} 
For simplicity, we set
\[\tilde{T} := W + \mu(T \otimes H^0(K_Y(3)))\subset H^0(K_Y \otimes L(3)).\]
Notice that the multiplication map
\[\mu : \tilde{T} \otimes H^0(L(-3)) \rightarrow H^0(K_Y \otimes L^2) \]
cannot be surjective, otherwise, as in the proof of Lemma \ref{omega2}, we get that 
\[ \pi_2 \circ \mu(\tilde{T} \otimes H^0(L(-3))) = H^{1,1}_{\rm ev}(X,\mathbb{C}) \]
and, given the first property of $W$, the latter equality implies the contradiction
\[\pi_2 \circ \mu(T \otimes H^0(K_Y \otimes L))) =
H^{1,1}_{\rm ev}(X,\mathbb{C}).\]
Now, by Proposition \ref{el}, we get that 
\[{\rm codim} \ \mu(T \otimes H^0(K_Y(3))) \geq \begin{cases} d - 5 + \alpha_Y - \beta_Y & {\rm
if \ L \ is \ (-d)-regular} \\ d - 5 & {\rm if} \ L = K_Y + dH + A \end{cases}. \]
$\hfill\square$
\newline
Therefore it will be enough to devise a mechanism for extracting codimension bounds for $T$ from
codimension bounds for $\mu(T\otimes H^0(K_Y(3)))$. This is the subject of the next section.
\newline 
We end the section by studying the vanishing of $H^1(\Omega_Y^2 \otimes L)$.
\begin{remark} If $d \geq 3\beta_Y - 3\alpha_Y + 13$ and $L$ is (-d)-regular or if $d \geq 2b_Y -
2a_Y + 13$ and $L = K_Y + dH + A$, then $H^1(\Omega_Y^2 \otimes L) = 0$.
\newline  
{\it Proof.} {\rm We just apply Griffiths' vanishing theorem \cite{gr} to the globally generated
vector bundle $E = \Omega_Y^2(3)$. We write
\[\Omega_Y^2 \otimes L = E ({\rm det}E + K_Y + B)\]
whence we just need to prove that
\[B = L - 12 H - 3 K_Y \]
is ample. By definition of $a_Y, b_Y, \alpha_Y$ and $\beta_Y$ we can write
\[-K_Y = (a - b)H + A'\]
where $A'$ is nef and $a = \alpha_Y, b = \beta_Y$ if $L$ is (-d)-regular, while
$a = a_Y, b = b_Y$ if $L = K_Y + dH + A$. Hence
\[B = (d - 12 - ub + u a)H + A''\]
where $A''$ is nef and $u = 2$ if $L = K_Y + dH + A$, $u = 3$ if $L$ is (-d)-regular. Therefore
$B$ is ample.} $\hfill\square$
\end{remark}
\begin{remark}
\label{quadric}
{\rm Notice that if $Y$ is a quadric hypersurface in $\mathbb{P}^4$, since $K_Y = -3H$, if $L =
(d-3)H$, we have that $H^1(\Omega_Y^2 \otimes L) = 0$ for $d \geq 7$, whence}
\[{\rm codim} \ T \geq d - 5.\]
\end{remark}
\section{Macaulay-Gotzmann for CM regular sheaves.}
We start by reviewing the situation for $\mathbb{P}^n$, which we will then generalise to 
arbitrary varieties.
\newline\newline
{\bf Definition of $c^{<d>}$ and $c_{<d>}$.} Given integers $c \geq 1, d \geq
1$, there  exists a unique sequence of integers $k_d, k_{d-1}, \ldots, k_f$ with $d \geq f \geq
1$  ($f$ is uniquely determined by $c$ and $d$) such that
\begin{enumerate}
\item $k_d > k_{d-1} > \ldots > k_f \geq f$,
\item $c = \sum\limits_{i = d}^f {k_i \choose i}$.
\end{enumerate}
Here and below we use the convention ${m \choose p} =  0$ if $m < p$.
We define
\[c^{<d>} := \sum_{i = d}^f {k_i + 1 \choose i + 1}\]
\[c_{<d>} := \sum_{i = d}^f {k_i - 1 \choose i}.\]
When $c = 0$ we set $c^{<d>} = c_{<d>} = 0$.
\newline
We have the following result of Macaulay and Gotzmann, which can be found in
\cite{gotz}, pages 64-65.
\newline\newline
{\bf Theorem (Macaulay, Gotzmann)}
{\it Let $V \subset H^0(\mathcal{O}_{\mathbb{P}^n}(d))$ be a subspace of codimension $c$. Then
the subspace
\[\mu(V\otimes H^0(\mathcal{O}_{\mathbb{P}^n}(1))) \subset
H^0(\mathcal{O}_{\mathbb{P}^n}(d + 1)) \] 
is of codimension at most $c^{<d>}$.}
\newline\newline 
Gotzmann proved the Macaulay-Gotzmann inequality using combinatorial algebraic techniques. 
Green gave a geometric proof in \cite{gmg}. We will now generalise the argument given by 
Green in order to prove that the Macaulay-Gotzmann inequality is valid for arbitrary 
Castelnuovo-Mumford regular sheaves.
\begin{result}
\label{cm}
Let $M$ be a Castelnuovo-Mumford regular coherent sheaf on a projective space $\mathbb{P}^N$. 
For $d \geq 1$ let 
\[V \subset H^0(M(d))\] 
be a subspace of codimension $c$, and define $V^{d+1} \subset H^0(M(d+1))$ by
\[V^{d+1} = \mu(V \otimes H^0(\mathcal{O}_{\mathbb{P}^N}(1))).\] 
Then 
\[{\rm codim} \ V^{d+1} \leq c^{<d>}.\]
\end{result}
The Theorem will follow from the following proposition.
\begin{proposition}
\label{onH}
Suppose that $V$, $M$ and $d$ are as above. Let $H$ be a generic hyperplane of $\mathbb{P}^N$
and denote by $M_H$ the restriction of $M$ to $H$. We further denote the restriction of $V$ to 
$H^0(M_H(d))$ by $V_H$. Then 
\[{\rm codim} \ V_H \leq c_{<d>}.\]
\end{proposition}
{\bf Proof.}
We shall proceed by a double induction on the dimension of the support of $M$ and
the number $d$. We assume now that $d \geq 2, {\rm dim Supp}(M) \geq 1$. The proof of the
Proposition for $d = 1$ or for sheaves with zero-dimensional supports is to be found in
subsections \ref{d=1} and  \ref{dimsupp0}.
\newline\newline
Let $H$ and $H'$ be two generic hyperplanes. We define the spaces $V^H$ (respectively $V^{H'}$) 
in the following way. Let $L_H$ (resp. $L_{H'}$) be a linear polynomial defining $H$ (resp.
$H'$). We define
\[V^H \subset H^0(M(d-1))\] by 
\[ v \in V^H \Leftrightarrow L_H \times v \in V.\]
(Similarly, $V^{H'}$ is defined by $v \in V^{H'} \Leftrightarrow L_{H'} \times v \in V$.)
We now consider the following exact sequence
\[0 \rightarrow H^0(M(d-1)) \stackrel{\times L_H}{\rightarrow }
H^0(M(d)) \stackrel{\rm res}{\rightarrow} H^0(M_H(d))\rightarrow 0.\]
Here, of course, we have right exactness of the sequence only because $M$ is a 
Castelnuovo-Mumford regular sheaf. There is an induced exact sequence
\[0 \rightarrow V^H \rightarrow V \rightarrow V_H \rightarrow 0 \]
whence we see that
\[{\rm codim} \ V = {\rm codim} \ V^H + {\rm codim} \ V_H.\]
We now consider the following commutative diagram
$$\xymatrix{ &  & 0 \ar[d] & 0 \ar[d] & \\
0 \ar[r] & (V^{H'})^{H} \ar[r] & V^{H'} \ar[r] \ar[d]^{\times L_{H'}} & 
(V^{H'})_{H} \ar[r] \ar[d]^{\times L_{H \cap H'}} & 0 \\
0 \ar[r] & V^H \ar[r] & V \ar[r]\ar[d]^{\rm res} & V_H
\ar[r]\ar[d]^{\rm res}& 0\\
0 \ar[r] & (V_{H'})^{H\cap H'}\ar[r] & V_{H'}\ar[r]\ar[d] & 
(V_{H'})_{H\cap H'}\ar[d]
\ar[r] & 0 \\ &  & 0 & 0 & \\}$$
In the above diagram, all the rows are exact (since $M_H$ is Castelnuovo-Mumford regular on 
$H$), as is the middle column. It is not immediate that the right-hand column is exact, but we will
be able to show that it is close enough to exact for our purposes.
\newline\newline
More precisely, 
\[(V_{H'})_{H \cap H'} = V_{|_{H \cap H'}} = (V_H)_{H \cap H'}\]
and hence the restriction map $V_H \rightarrow (V_{H'})_{H \cap H'}$ is a
surjection. We have automatically that
\[(V^{H'})_H \subset (V_H)^{H\cap H'}\]
and hence the composition of the maps $\times L_{H \cap H'}$ and ${\rm res}$ is zero.
It follows that
\[{\rm codim} \ V_H \leq {\rm codim} \ (V_{H'})_{H \cap H'} + {\rm codim} \ (V^{H'})_H.\]
We denote by $c'$ the codimension of $V_H$ for generic $H$. Hence, since $H'$ has been chosen 
generic, ${\rm codim} \ V_{H'} = c'$. We have that 
\[{\rm codim} \ V^{H'} = c - c'.\] 
We note that
\begin{enumerate}
\item $V^{H'} \subset H^0(M(d-1))$ and hence by the induction hypothesis 
\[{\rm codim} \ (V^{H'})_H \leq (c - c')_{<d - 1>}.\]
\item The dimension of the support of $M_{H'}$ is strictly less than the dimension of the 
support of $M$ and hence by the induction hypothesis
\[{\rm codim} \ (V_{H'})_{H \cap H'} \leq c'_{<d>}.\]
\end{enumerate}
It follows that
\[c'\leq c'_{<d>}+ (c-c')_{<d-1>}.\]
Green shows in \cite{gmg}, pages 77-78, that this inequality implies that $c'\leq c_{<d>}$.
\newline\newline
It remains only to prove the Proposition for zero-dimensional sheaves or for $d=1$. 
\subsubsection{The case d=1.}
\label{d=1}
We have that for any $c \neq 0$, 
\[c_{<1>} = c-1.\] 
We suppose first that $V \neq H^0(M(1))$. If for generic $H$ 
\[{\rm codim} \ V_H > c_{<1>}\] 
then for generic $H$ 
\[V^H = H^0(M).\] 
In other words, for generic $H$
\[L_H \times H^0(M) \subset V.\]
It follows that
\[\mu(H^0(M), H^0(\mathcal{O}_{\mathbb{P}^N}(1))) \subset V.\] 
Since $M$ is Castelnuovo-Mumford regular, it follows that $V = H^0(M(1))$ which contradicts our 
supposition that  $V \neq H^0(M(1))$. 
\newline\newline
But if $c=0$ then $c_{<1>} = 0$ and Proposition \ref{onH} is immediate. This completes the
proof of the Proposition in the case where $d=1$.
\subsubsection{The case where the dimension of the support of $M$ is zero.}
\label{dimsupp0}
In this case, for generic $H$, $H^0(M_H(d))=0$, and hence ${\rm codim}\ V_H = 0$. This completes 
the proof of the Proposition in the case where the dimension of the support of $M$ is zero.
\newline\newline
This completes the proof of Proposition \ref{onH}. $\hfill\square$
\newline\newline
We now show how Proposition \ref{onH} implies Theorem \ref{cm}. We proceed by induction on the
dimension of  the support of $M$. We consider the following exact sequence, where $H$ is once
again a generic  hyperplane in $\mathbb{P}^N$,
\[ 0\rightarrow (V^{d+1})^H \rightarrow V^{d+1} \rightarrow (V^{d+1})_H \rightarrow 0\]
from which it follows that
\[{\rm codim} \ V^{d+1} = {\rm codim} \ (V^{d+1})^H + {\rm codim} \ (V^{d+1})_H.\]
We note that $V \subset (V^{d+1})^H$ and $(V_H)^{d+1} \subset (V^{d+1})_H$ from which it follows
that
\[{\rm codim} \ V^{d+1} \leq c + (c_{<d>})^{<d>} \leq c^{<d>}.\]
This completes the proof of Theorem \ref{cm}. $\hfill\square$
\section{Proof of the main theorems.}
We will now show how all this ties together to give a proof of the main theorems.
We henceforth set
\[a = \begin{cases} \alpha_Y & {\rm if \ L \ is \ (-d)-regular} \\ a_Y & {\rm if} \ L = K_Y +
dH + A \end{cases}, \ b = \begin{cases} \beta_Y & {\rm if \ L \ is \ (-d)-regular} \\ b_Y & {\rm
if} \ L = K_Y + dH + A \end{cases}\]
where $\alpha_Y, \beta_Y, a_Y$ and $b_Y$ are as in Definitions \ref{alpha} and \ref{a}.
\newline
It is now that we will use the supposition that $(Y, H)$ is not a linear $\mathbb{P}^2$-bundle,
hence $K_Y(3)$ is very ample, or, alternatively, that $a \leq 3$ (the case of the quadric is done by
Remark \ref{quadric}). The case $a = 4$ will be dealt with at the end of the article.
\newline
We start with the following lemma.
\begin{little}
\label{T'}
Suppose $d \geq 5$ and let $T \subset H^0(L)$ be of codimension $c \leq d - 4$. Define 
\[T' := \mu (T \otimes H^0(\mathcal{O}_Y(3-a))) \subset H^0(L(3-a)).\]
Then 
\[{\rm codim} \ T' \leq c\]
\end{little}
{\bf Proof.}
When $L$ is (-d)-regular we can write 
\[L = M(d),\]
where $M$ is a Castelnuovo-Mumford regular sheaf. Also when $L = K_Y + dH + A$, since
$M := K_Y + 4H + A$ is Castelnuovo-Mumford regular, we can write 
\[L = M(d - 4),\] 
where $M$ is a Castelnuovo-Mumford regular sheaf. Applying $(3-a)$-times Theorem \ref{cm},
we obtain the result. $\hfill\square$
\newline\newline
We denote now by $n$ the integer $\lfloor \frac{d+3-a}{b}\rfloor - 4$. 
We will also denote the very ample line bundle $K_Y(a)$ by $P$, and the bundle $L(3-a)$ by 
$L'$. We have the following lemma.
\begin{little}
\label{L'}
The line bundle $L'$ can be written in the form
\[L' = M_P + nP \]
where $M_P$ is a sheaf which is Castelnuovo-Mumford regular {\it with respect to the projective 
embedding defined by P.}
\end{little}
{\bf Proof.}
We know by definition of $a$ and $b$ that there is a nef line bundle $N$ such that
\[bH = K_Y + aH + N\]
from which it follows that
\[(d + 3 - a)H = (n + 4)P + (n + 4)N + rH\]
for some $r \geq 0$, hence
\[(d + 3 - a)H = (n + 4)P + A'\]
where  $A'$ is a nef line bundle. Now
\[M_P := L' - n P =  \begin{cases} 4P + A_1 & {\rm if \ L \ is \ (-d)-regular} \\ 
K_Y + 4 P + A_2 & {\rm if} \ L = K_Y + dH + A \end{cases}\]
for some nef line bundles $A_1, A_2$. This clearly implies, by Kodaira vanishing, that $M_P$ is
Castelnuovo-Mumford regular with respect to $P$ in the case $L = K_Y + dH + A$. But also in the
other case, for each $1 \leq i \leq 3$, we can write
\[M_P - iP = K_Y + aH + (3 - i)P + A_1\]
whence again we have Castelnuovo-Mumford regularity by Kodaira vanishing since now 
$a = \alpha_Y > 0$ by definition. $\hfill\square$
\newline\newline
We are now in a position to prove the following proposition.
\begin{proposition} 
\label{overlineT}
Suppose $d \geq 5$ and let $T \subset H^0(L)$ be of codimension $c \leq d - 4$. Define 
\[\overline{T} := \mu (T \otimes H^0(K_Y(3)) \subset H^0(K_Y \otimes L(3)). \]
Then 
\[{\rm codim} \ \overline{T} \leq c^{<n>}.\]
\end{proposition}
{\bf Proof.}
With $T'$ as in Lemma \ref{T'}, we note that
\[\mu(T' \otimes H^0(K_Y(a)) \subset \overline{T}.\]
We know by Lemma \ref{T'} that ${\rm codim} \ T' \leq c$. We know further by Lemma \ref{L'}
that
\[L' = M_P + nP \]
and hence Theorem \ref{cm} applied to the map
\[\mu : T' \otimes H^0(P) \rightarrow H^0(K_Y \otimes L(3))\]
gives us that
\[{\rm codim} \ \mu(T' \otimes H^0(K_Y(a)) \leq c^{<n>}.\]
From this it follows that 
\[\hskip 4cm {\rm codim} \ \overline{T} \leq c^{<n>}. \hskip 5.8cm \square \]
\newline
By Proposition \ref{firstlower} we know that 
\[{\rm codim} \ \overline{T} \geq \begin{cases} d - 5 + \alpha_Y - \beta_Y & {\rm if \ L \ is \
(-d)-regular} \\ d - 5 & {\rm if} \ L = K_Y + dH + A \end{cases} \] 
and hence either $c \geq d - 3$ or
\[c^{<n>} > \begin{cases} d - 6 + \alpha_Y - \beta_Y & {\rm if \ L \ is \
(-d)-regular} \\ d - 6 & {\rm if} \ L = K_Y + dH + A \end{cases}. \]
The following elementary lemma will allow us to control the growth of $c^{<n>}$.
\begin{little}
\label{growth}
If there exists an integer $e \geq 0$ such that
\[c < \sum_{i=0}^e (n + 1 - i) \] 
then $c^{<n>} \leq c + e$.
\end{little}
{\bf Proof.}
The Lemma being obvious for $c = 0$ we suppose $c \geq 1$ and 
$c = \sum\limits_{i = n}^f {k_i \choose i}$. Observe that
\[\sum_{i=0}^e (n + 1 - i) \leq \frac{(n + 1)(n + 2)}{2}.\]  
Now suppose $k_i = i$ for $f \leq i \leq f_1$ for some $f - 1 \leq f_1 \leq n$, $k_i = i +
1$  for $f_1 + 1 \leq i \leq f_2$ for some $f_2$ such that $f_1 \leq f_2 \leq n$ and
$k_i \geq i + 2$ for $f_2 + 1 \leq i \leq n$ (the case $f - 1 = f_1$ simply means that no $k_i$ is
equal to $i$, and similarly for $f_2$). Then, if
$f_2 < n$, we have
\[c \geq {k_n \choose n} \geq {n + 2 \choose 2} = \frac{(n + 1)(n + 2)}{2}\]
contradicting the hypothesis. Therefore $f_2 = n$ and
\[c^{<n>} = c + n - f_1\]
and it remains to show that $n - f_1 \leq e$. Since we can write
\[c = \sum_{i = 0}^{n - f_1} (n + 1 - i) - f \]
if $n - f_1 \geq e + 1$ we deduce the contradiction
\[\hskip 4cm c \geq \sum_{i=0}^e (n + 1 - i). \hskip 4.8cm \square \] 
\newline\newline
In particular, it follows that
\begin{little}
Suppose $L = K_Y + dH + A$, $b_Y \geq 2$ and
\[d - 6 - b_Y < \sum_{i = 0}^{b_Y} (n + 1 - i).\] 
Then 
\[{\rm codim} \ T >  d - 6 - b_Y.\]
If $b_Y = 1$, then 
\[{\rm codim} \ T >  d - 6.\]
\end{little}
{\bf Proof.}
By Lemma \ref{growth}, if $b_Y \geq 2$,
\[ d - 6 - b_Y < \sum_{i = 0}^{b_Y} (n + 1 - i)\] 
and  
\[c = {\rm codim} \ T \leq d - 6 - b_Y\] 
then, by Proposition \ref{overlineT},
\[{\rm codim} \ \overline{T} \leq c^{<n>} \leq d  - 6.\]
But this is impossible by Proposition \ref{firstlower}. If $b_Y = 1$ and $c \leq d - 6$
we have $c \leq n$ hence 
\[{\rm codim} \ \overline{T} \leq c^{<n>} = c \leq d  - 6,\]
again impossible by Proposition \ref{firstlower}. $\hfill\square$
\newline\newline
Similarly we have
\begin{little}
Suppose $L$ is (-d)-regular, $\beta_Y \geq 2$ and
\[d - 6 + \alpha_Y - 2 \beta_Y < \sum_{i = 0}^{\beta_Y} (n + 1 - i).\] 
Then 
\[{\rm codim} \ T >  d - 6 + \alpha_Y - 2\beta_Y.\]
If $\beta_Y = 1$, then 
\[{\rm codim} \ T >  d - 7 + \alpha_Y.\]
\end{little}
We now require only the following lemma.
\begin{little}
If $b_Y \geq 2$ and $d \geq \frac{b_Y(b_Y^2 + 7 b_Y - 6)}{2}$ then 
\begin{equation}
\label{eq} 
d - 6 - b_Y < \sum_{i = 0}^{b_Y} (n + 1 - i).\end{equation} 
If $\beta_Y \geq 2$ and $d \geq \frac{\beta_Y^2 (\beta_Y + 5)}{2}$ then
\[d - 6 + \alpha_Y - 2 \beta_Y < \sum_{i = 0}^{\beta_Y} (n + 1 - i).\] 
\end{little}
{\bf Proof.}
We note first that $n \geq \lfloor\frac{d}{b_Y} \rfloor - 4$ and it
follows that 
\[b_Y(n + 1) > d - 4 b_Y.\] 
Hence we have that
\[\sum_{i = 0}^{b_Y} (n + 1 - i) > d - 4 b_Y + (n + 1) - \frac{b_Y(b_Y + 1)}{2}.\]
In particular, if 
\[d - 6 - b_Y \leq d - 4 b_Y + (n + 1) - \frac{b_Y(b_Y + 1)}{2}\]
then (\ref{eq}) is immediately satisfied.
This inequality is equivalent to
\[- 7 + 3b_Y \leq n - \frac{b_Y(b_Y + 1)}{2}\]
and since $n \geq \lfloor \frac{d}{b_Y} \rfloor - 4$, (\ref{eq}) will be
satisfied provided that
\[- 7 + 3 b_Y \leq \lfloor \frac{d}{b_Y} \rfloor - 4 - \frac{b_Y(b_Y + 1)}{2} \]
which is equivalent to
\[- 3 + 3 b_Y + \frac{b_Y(b_Y + 1)}{2} \leq \lfloor \frac{d}{b_Y} \rfloor\] 
which is equivalent to
\[\frac{b_Y(b_Y^2 + 7 b_Y - 6)}{2} \leq d.\]
The second assertion of the Lemma is proved similarly. $\hfill\square$
\newline\newline
{\bf Completion of the proof of Theorems \ref{(-d)-regular} and \ref{pluriadjoint}.}
\newline
The results proved so far (together with Remark \ref{quadric}) give a proof of the Theorems under
the hypothesis that $(Y, H)$ is not a linear $\mathbb{P}^2$-bundle. In the latter case since
$K_Y(4)$ is very ample, repeating verbatim the whole proof replacing everywhere $K_Y(3)$ with
$K_Y(4)$ and using $a_Y = \alpha_Y = 4$ we get the desired bound. $\hfill\square$
\newline\newline
{\bf Proof of Corollary \ref{blow-up}.}
\newline
This is a straightforward generalisation of \cite{an} given the
following two facts : 
\newline
1. a lower bound on the codimension on the components of the Noether-Lefschetz locus ${\rm
NL}(\mathcal{O}_Y(d))$ that insures that they have codimension at least two (our hypothesis
$d \geq 7 + e$); 
\newline
2. the fact that, on a general surface $X$ not in ${\rm NL}(\mathcal{O}_Y(d))$ we have that if a
complete intersection of $X$ with another surface in $|\mathcal{O}_Y(d)|$ is reducible then its
irreducible components are also complete intersection of $X$ with another surface in
$|\mathcal{O}_Y(s)|$ for some $s$ (this is needed in the proof of \cite[Prop.\ 2.1]{an} and is insured,
in our case, by the hypothesis ${\rm Pic}(Y) \cong \mathbb{Z}$). $\hfill\square$

{\bf Addresses of the authors:}
\newline\newline
Angelo Felice Lopez, Dipartimento di Matematica, Universit\`a di Roma Tre,
Largo San Leonardo Murialdo 1, 00146, Roma, Italy.\newline
e-mail {\tt lopez@mat.uniroma3.it}
\newline\newline
Catriona Maclean, UFR de Math\`ematiques, UMR 5582 CNRS / Universit\`e J. Fourier
100, rue des Maths, BP74, 38402 St Martin d'Heres, France.\newline
e-mail {\tt Catriona.MacLean@ujf-grenoble.fr}

\end{document}